\documentclass[10pt,a4paper]{proc-l}
\usepackage{amsmath}
\usepackage{latexsym}
\usepackage{amsfonts}
\usepackage{amssymb}
\usepackage{amsthm,amsxtra}
\usepackage[all,web,colour]{xy}

\newtheorem{theorem}{Theorem}

\newtheorem{proposition}[theorem]{Proposition}
\newtheorem{corollary}[theorem]{Corollary}

\theoremstyle{definition}

\theoremstyle{remark}
\newtheorem{remark}{Remark}

\newcommand{\Z}{\mathbb Z}

\newcommand{\C}{\mathbb C}
\newcommand{\+}{\!+\!}

\newcommand{\PV}{${\rm P}_{\rm V}\;$}
\newcommand{\PVI}{${\rm P}_{\rm VI}\;$}
\newcommand{\PaII}{${\rm P}_{\rm II}\;$}
\newcommand{\PIII}{${\rm P}_{\rm III}\;$}

\newcommand{\PIIIprime}{${\rm P}_{\rm III^{\prime}}\;$}
\newcommand{\threehalf}{
        {\lower0.00ex\hbox{\raise.6ex\hbox{\the\scriptfont0 3}
                           \kern-.5em\slash\kern-.1em\lower.45ex
                                     \hbox{\the\scriptfont0 2}}}}
\newcommand{\half}{
        {\lower0.00ex\hbox{\raise.6ex\hbox{\the\scriptfont0 1}
                           \kern-.5em\slash\kern-.1em\lower.45ex
                                     \hbox{\the\scriptfont0 2}}}}
\newcommand{\quarter}{
        {\lower0.00ex\hbox{\raise.6ex\hbox{\the\scriptfont0 1}
                           \kern-.5em\slash\kern-.1em\lower.45ex
                                     \hbox{\the\scriptfont0 4}}}}
\newcommand{\eighth}{
        {\lower0.00ex\hbox{\raise.6ex\hbox{\the\scriptfont0 1}
                           \kern-.5em\slash\kern-.1em\lower.45ex
                                     \hbox{\the\scriptfont0 8}}}}
\newcommand{\sixteenth}{
        {\lower0.00ex\hbox{\raise.6ex\hbox{\the\scriptfont0 1}
                           \kern-.5em\slash\kern-.1em\lower.45ex
                                     \hbox{\the\scriptfont0 16}}}}
\newcommand{\thirtytwo}{
        {\lower0.00ex\hbox{\raise.6ex\hbox{\the\scriptfont0 1}
                           \kern-.5em\slash\kern-.1em\lower.45ex
                                     \hbox{\the\scriptfont0 32}}}}

\begin{document}
\title{New Transformations for Painlev\'e's Third Transcendent}

\author{N.S.~Witte}
\address{Department of Mathematics and Statistics,
University of Melbourne, Victoria 3010, Australia}
\email{N.Witte@ms.unimelb.edu.au}

\date{January 26, 2002 and, in revised form, XXX XX, 2002.}

\begin{abstract}
We present transformations relating the third transcendent of Painlev\'e 
with parameter sets located at the corners of the Weyl chamber for the
symmetry group of the system, the affine Weyl group of the root system 
$ B^{(1)}_2 $, to those at the origin. This transformation entails a scaling
of the independent variable, and implies additive identities for the 
canonical Hamiltonians and product identities for the $\tau$-functions
with these parameter sets.
\end{abstract}

\commby{Mark J. Ablowitz}
\maketitle

A curious anomaly has existed in the theory of the Painlev\'e transcendents
for many years, for almost a century in fact. This anomaly is the existence of
a particular transformation for the second Painlev\'e transcendent 
$ q(t;\alpha +\half) $ which has no known analogue for any of the other 
transcendents.  One recalls the second Painlev\'e transcendent 
$ q=q(t;\alpha +\half) $
is by definition a solution of the \PaII differential equation
\begin{equation}
   {d^2q \over dt^2} = 2q^3+tq+\alpha\ ,
\end{equation}
with parameter $ \alpha \in {\mathbb C}$.
The transformation, found by Gambier in 1909 \cite{Ga_1909}, relates solutions 
of the \PaII transcendent with parameter values at
$ v_1 := \alpha+\half = 0 $ and $ v_1 = \half $ (also between the solutions at 
$ v_1 = \half $ and $ v_1 = 1 $ but this can be constructed from the previous
one and the B\"acklund transformation $ v_1 \mapsto v_1\+ 1 $).
The B\"acklund symmetries of \PaII are a realisation of the affine Weyl 
group associated with the root system $ A^{(1)}_1 $, and using 
\begin{equation}
\begin{split}
   v_1 & \mapsto v_1 + 1, \quad\text{shift} \cr
   v_1 & \mapsto -v_1, \quad\text{reflection},
\end{split}
\label{PII_Bxfm}
\end{equation}
transformations of this group,
solutions for arbitrary $ v_1 $ can be related to a solution whose parameter
$ v_1 \in (0,\half] $ (or in $ [\half,1) $). This interval defines the Weyl
chamber and the above transformation in question relates the solutions at the 
endpoints of this interval (see Fig. 1). 
%

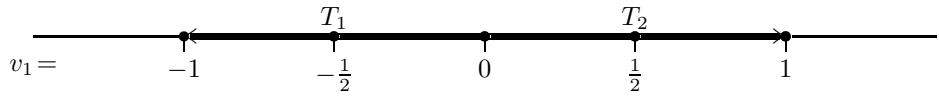
\begin{figure}[H]
\begin{center}

  \newcommand{\monad}{
   \save 
  {0*{\bullet};<2cm,0cm>**@{-}
  }\restore}

\[
 \renewcommand{\latticebody}{\drop{\monad}}
 \begin{xy}
    *\xybox{0;<2.0cm,0cm>:<1.0cm,1.73205081cm>::
           ,0,{\xylattice{-1}{2}{0}{0}},
	   ,{\ar@{=>}_{\displaystyle T_1} 0;0+(-2,0)}
	   ,{\ar@{=>}^{\displaystyle T_2} 0;0+(2,0)}
           ,0+(-3,0)="O";"O"+<2.0cm,0cm>**@{-}
 	   ,"O";"O"+< .0cm,-.2cm>+D*+!U{v_1\!=}
           ,0+(-2,0)="O"*{\bullet};"O"+<2.0cm,0cm>**@{-}
 	   ,"O";"O"+< .0cm,-.2cm>**@{-}+D*+!U{-1}
           ,0+(-1,0)="O"*{\bullet};"O"+<2.0cm,0cm>**@{-}
 	   ,"O";"O"+< .0cm,-.2cm>**@{-}+D*+!U{-\frac{1}{2}}
           ,0+( 0,0)="O"*{\bullet};"O"+<2.0cm,0cm>**@{-}
 	   ,"O";"O"+< .0cm,-.2cm>**@{-}+D*+!U{0}
           ,0+( 1,0)="O"*{\bullet};"O"+<2.0cm,0cm>**@{-}
 	   ,"O";"O"+< .0cm,-.2cm>**@{-}+D*+!U{\frac{1}{2}}
           ,0+( 2,0)="O"*{\bullet};"O"+<2.0cm,0cm>**@{-}
 	   ,"O";"O"+< .0cm,-.2cm>**@{-}+D*+!U{1}
           }
 \end{xy} \]

\caption{Parameter space $ (v_1) $ associated
with the simple roots of the root system $ A_{1}^{(1)} $.}
\label{A1-fig}
\end{center}
\end{figure}


%
Specifically the transformation takes the form
\begin{equation}
\begin{split}
  T & = -2^{-1/3}t
    \\
  2^{1/3}Q^2(T) & = {dq(t) \over dt}+q^2(t)+\half t \ ,
\end{split}
\label{PII_xfm}
\end{equation}
where $ q(t) := q(t;0),\; Q(T) := q(T;\half) $ and its inverse is
\begin{equation}
\begin{split}
  t & = -2^{1/3}T
    \\
  q(t) & = -2^{-1/3}{1 \over Q(T)}{dQ(T) \over dT} \ .
\end{split}
\label{PII_ixfm}
\end{equation}
We draw the reader's attention to some pertinent observations that distinguish
this transformation from the B\"acklund (or equivalently Schlesinger) 
transformations (\ref{PII_Bxfm}) -- 
firstly the shift in the parameters is always the smallest possible shift, one half
unit, on 
the root lattice whereas the B\"acklund transformation shift is twice this, 
secondly the independent variable undergoes a scale change in contrast to 
remaining a fixed variable, and thirdly the transformation is quadratic in
$ q $ and $ Q $ but linear in the derivatives.
The purpose of this note is to demonstrate that such a transformation exists for 
\PIIIprime as well, and to illustrate the many analogous properties that the
new transformation for \PIIIprime has with those of the \PaII transformation. 
After recognising these similarities one may well talk of a class of such
transformations for all systems, \PaII to \PVI.
Actually an analogue of this class of transformation was noted by Okamoto 
\cite{Ok_1986} in the case of the fourth transcendent but appeared not to be 
appreciated as such.
Aside from the intrinsic interest of these transformations there are some
practical motivations from the theory of random matrices where they
arise as multiplicative identities for the probabilities that certain spectral
intervals are free of eigenvalues \cite{Wi_2001a}.

We first present the context for these transformations by reviewing the 
understanding of the B\"acklund symmetries of \PIII and the underlying affine 
Weyl group structure.
The \PIIIprime differential equation for $ q(t) $, in the convention of Okamoto
\cite{Ok_1987c}, is
\begin{equation}
      {d^2 q \over dt^2} = {1 \over q} \Big ( {dq \over dt} \Big )^2 -
  {1 \over t} {dq \over dt} + {q^2 \over 4t^2} (\gamma q + \alpha) +
  {\beta \over 4t}  + {\delta \over 4q} \ ,
\label{PIII_ode}
\end{equation}
with the parameter identifications
\begin{equation}
  \alpha = -4\eta_{\infty}v_2,\quad \beta = 4\eta_{0}(v_1+1),\quad
  \gamma = 4\eta^2_{\infty},\quad \delta = -4\eta^2_{0} \ .
\end{equation}
Here $ \eta_{\infty}, \eta_{0} $ are fixed at unity and all discussions of 
parameters will be in terms of the coordinates $ {\bf v}=(v_1,v_2) \in \C^2 $ on the 
root lattice.
There is a Hamiltonian structure underlying the \PIIIprime system 
$ \{q,p;t,H\} $ and the Hamiltonian is conventionally taken to be
\begin{equation}
  tH = q^2 p^2 - (\eta_{\infty}q^2 + v_1q - \eta_{0}t)p
               + \half \eta_{\infty}(v_1 + v_2) q\ .
\label{PIII_H}
\end{equation}
The Hamilton equations of motions are then
\begin{align}
     tq' & = 2q^2p - \eta_{\infty}q^2 - v_1 q + \eta_{0}t
         \\
     tp' & = -2qp^2 + (2\eta_{\infty}q + v_1)p - \half \eta_{\infty}(v_1 + v_2)
     \ ,
\label{PIII_Heom}
\end{align}
and $ q(t) $ satisfies (\ref{PIII_ode}).
Okamoto \cite{Ok_1987c} and Kajiwara et al \cite{KMNOY_2001} have identified two 
B\"acklund or Schlesinger transformations with the action on the parameters
\begin{equation}\label{T}
  T_1\cdot\mathbf{v} = (v_1+1, v_2+1), \qquad
  T_2\cdot\mathbf{v} = (v_1+1, v_2 - 1).
\end{equation}
The operators $T_1$ and $T_2$ can be constructed out of the generators
$s_0, s_1, s_2$ which are reflections associated with the underlying 
$ B_2 $ root lattice, whose action (following \cite{Ok_1987c} and \cite{KMNOY_2001}) 
on $ \mathbf{v}, q,p $ is given in Table \ref{t2}. According to Table \ref{t2} we 
have
\begin{equation}
  T_1 = s_0 s_2 s_1 s_2, \qquad
  T_2 = s_2s_0s_2s_1.
\label{PIII_shifts}
\end{equation}
See Fig. 2 for the $ B_2 $ root lattice.

\begin{table}
\begin{center}
\begin{tabular}{|c||c|c|c|c|c|}\hline
 & $v_1$ & $v_2$ & $p$ & $q$ & $t$  \\ \hline
 $s_0$ & $-1-v_2$ & $-1-v_1$  &
 ${\displaystyle q \over \displaystyle t}
   \left[q(p-1) - {1 \over 2}(v_1-v_2)\right] + 1$ &
 $-{\displaystyle t \over \displaystyle q} $ & $t$ \\
 $s_1$ & $v_2$ & $v_1$ &
 $p$ & $q + {\displaystyle v_2-v_1 \over \displaystyle 2(p-1)} $ & $t$ \\
 $s_2$ & $v_1$  & $-v_2$ &  $1-p$ & $-q$ &
 $-t$   \\ \hline
\end{tabular}
\end{center}
\caption{\label{t2}  Generators of the B\"acklund transformations realising a
representation of the extended affine Weyl group for the root system 
$ B^{(1)}_2 $.}
\end{table}
%


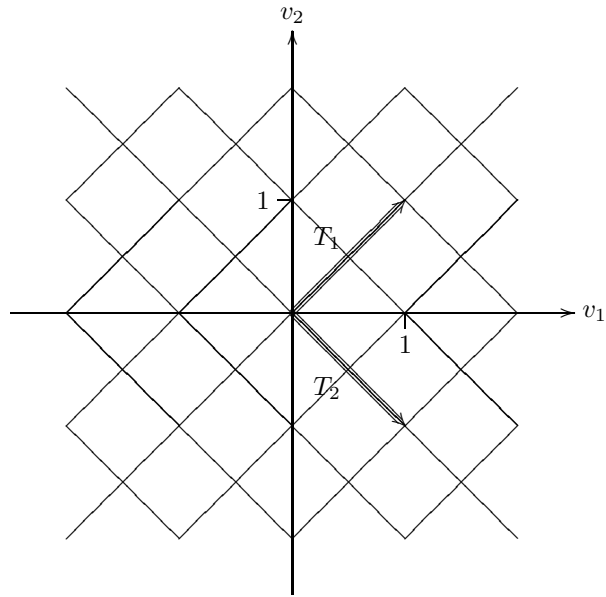
\begin{figure}[H]
\begin{center}

  \newcommand{\diad}{
   \save 
  { 0;<+1.5cm,+1.5cm>**@{-}
   ,0;<+1.5cm,-1.5cm>**@{-}
   ,0+<0.0cm,1.5cm>="O";"O"+<+1.5cm,+1.5cm>**@{-}
   ,"O";"O"+<+1.5cm,-1.5cm>**@{-}
  }\restore}

\[
 \renewcommand{\latticebody}{\drop{\diad}}
 \begin{xy}
    *\xybox{0;<1.5cm,0cm>:<0.0cm,1.5cm>::
           ,0,{\xylattice{-1}{2}{-0}{1}}="G",
 	   ,0+( 1,-1.5)="O",{\ar@{->} "O";"O"+(0, 5)},"O"+(0, 5)*+!D{v_2}
 	   ,0+(-1.5, 1)="O",{\ar@{->} "O";"O"+(5, 0)},"O"+(5, 0)*+!L{v_1}
 	   ,0+(-1, 3)="O";"O"+(1,-1)**@{-}
 	   ,0+( 0, 3)="O";"O"+(1,-1)**@{-}
 	   ,0+( 1, 3)="O";"O"+(1,-1)**@{-}
 	   ,0+( 2, 3)="O";"O"+(1,-1)**@{-}
 	   ,0+(-1,-1)="O";"O"+(1, 1)**@{-}
 	   ,0+( 0,-1)="O";"O"+(1, 1)**@{-}
 	   ,0+( 1,-1)="O";"O"+(1, 1)**@{-}
 	   ,0+( 2,-1)="O";"O"+(1, 1)**@{-}
           ,{\ar@{=>}^{\displaystyle T_1} 0+(1,1);0+(2,2)}
           ,{\ar@{=>}_{\displaystyle T_2} 0+(1,1);0+(2,0)}
           ,"O"+(0,2);"O"+(0,2)+< .0cm,-.2cm>**@{-}+D*+!U{1}
           ,"O"+(-1,3);"O"+(-1,3)+<-.2cm, .0cm>**@{-}+D*+!R{1}
           }
 \end{xy} \]

\caption{
Parameter space $ (v_1, v_2) $ associated with the root system $ B_2^{(1)} $.}
\label{B2-fig}

\end{center}
\end{figure}

%
%
%


%
Associated with the Hamiltonian is the $ \tau $-function defined by
\begin{equation}
   H(t) := {d \over dt}\log\tau \ .
\label{tau-fn}
\end{equation}
It has the property of being entire in $ t \in \C\backslash\{0,\infty\} $ and 
having only simple zeros.
One-parameter solutions, a particular type of classical solution, of \PIIIprime are 
possible only when the parameters satisfy \cite{Ok_1987c,Mu_1995,Gr_1999}
\begin{equation}
   v_2 \pm v_1 \in 2\Z \ .
\end{equation}
The corresponding $ \tau $-functions have the explicit Toeplitz determinant 
form involving Bessel functions, 
\begin{equation}
  \det[ I_{\nu+j-k}(\sqrt{t})+c K_{\nu+j-k}(\sqrt{t}) ]_{0 \leq j,k \leq n}\ ,
\end{equation}
$ c, \nu \in \C $ or
\begin{equation}
  \det[ J_{\nu+j-k}(\sqrt{t})+c Y_{\nu+j-k}(\sqrt{t}) ]_{0 \leq j,k \leq n}\ .
\end{equation}
Rational, or zero-parameter, solutions are possible only when \cite{Gr_1999}
\begin{equation}
   v_2 \pm v_1 \in 2\Z + 1 \ .
\end{equation}
These are constructed by recurrences derived using the operators
(\ref{PIII_shifts}) from a seed solution.
Additional studies of the B\"acklund transformations and exact solutions 
of the third Painlev\'e equation from various perspectives are given 
in \cite{Lu_1967,Gr_1973,Gr_1975a,Gr_1975b,Ai_1979,BK_1985,MCB_1997,MW_1998,UW_1998,GJP_2001,GJP_2001a}.

Let us now return to Gambier's transformation (\ref{PII_xfm},\ref{PII_ixfm}).
Expressed in canonical variables (\ref{PII_xfm}) takes the form
\begin{equation}
\begin{split}
   T & = -2^{-1/3}t
   \\
   Q & = 2^{-1/6}\sqrt{p}
   \\
   P & = -2^{-1/3}\left[ \half t - p + q\sqrt{2p} \right]\ ,
\end{split}
\label{PII_Hxfm}
\end{equation}
and the inverse transformation (\ref{PII_ixfm}) is
\begin{equation}
\begin{split}
   t & = -2^{1/3}T
   \\
   q & = -2^{-1/3} {P-Q^2-T/2 \over Q}
   \\
   p & = 2^{1/3} Q^2\ .
\end{split}
\end{equation}
The corresponding transformation for \PIIIprime is given in the following result.

\begin{proposition}
The canonical variables of the \PIIIprime Hamiltonian system $ \{q, p; t, h\} $ 
with $ {\bf v}=(0,0) $ are related to the \PIIIprime Hamiltonian system 
$ \{Q, P; T, H\} $ with $ {\bf v}=(-1,0) $ by the direct canonical
transformation
\begin{equation}
\begin{split}
   T & = \quarter t
   \\
   Q & = i {\sqrt{t} \over 2} \left( \sqrt{p} + \sqrt{p-1} \right)
   \\
   P & = {i \over \sqrt{t}} 
            { \half + q\sqrt{p}\sqrt{p-1} - i\sqrt{t}\sqrt{p} \over
              \sqrt{p} + \sqrt{p-1} }\ .
\end{split}
\label{PIII_Wxfm}
\end{equation}
The inverse of this transformation is
\begin{equation}
\begin{split}
   t & = 4T
   \\
   q & = -2TQ{ 2Q^2-4Q^2P-Q-2T \over (Q^2+T)(Q^2-T) }
   \\
   p & = -{ (Q^2-T)^2 \over 4TQ^2 } \ .
\end{split}
\label{PIII_Wixfm}
\end{equation}
\end{proposition}

\begin{proof}[First Proof]
One can directly verify that $ Q(T) $ satisfies the differential equation 
(\ref{PIII_ode}) with appropriate parameters by hand or otherwise.
\end{proof}

\begin{proof}[Second Proof]
One can adapt Gromak's method for \PaII \cite{Gr_1999} to this case in the following 
sequence of steps. In the first step one eliminates $ q(t) $ instead of $ p(t) $
from Hamilton's equations of motion (\ref{PIII_Heom}) and arrives at 
\begin{equation}
   p'' = \half\left( {1 \over p}+{1 \over p-1} \right)p'{}^2
         - {p' \over t} + {2 \over t}p(1-p) 
         + {1 \over 4t^2}{ (v_1+v_2)^2(2p-1)-4v_1v_2 p^2 \over
                           2p(1-p) } \ .
\label{PV_dual}
\end{equation}
Employing the transformation
\begin{equation}
   p = {y \over y-1}\ ,
\end{equation}
one finds $ y(t) $ satisfies the standard form for Painlev\'e's fifth transcendent
with the parameters
\begin{equation}
  \alpha = \eighth(v_1-v_2)^2 ,\quad \beta = -\eighth(v_1+v_2)^2 ,\quad
  \gamma = 2,\quad \delta = 0 \; .
\end{equation}
and in this sense this particular form of \PV is dual to the general \PIII
\cite{GF_2001,Gr_1999,Gr_1984,Gr_1975a}.
In the second step we need to transform this into a \PIII transcendent but the crucial 
point is that one should not use a contact transformation to effect this under the most
general conditions but employ a point transformation which will only work for
particular parameters. By setting $ v_1 = v_2 = 0 $ such a transformation is
\begin{equation}
   y = \left({ u+1 \over u-1 }\right)^2\ ,
\end{equation}
which can be applied to \PV when $ \alpha = \beta = 0 $ and general 
$ \gamma, \delta $ \cite{Pa_1902}. 
The resulting differential equation for $ u(t) $ is that of a
\PIII transcendent with $ \alpha = -1/2, \beta = 1/2, \gamma = \delta = 0 $, which 
is in a nonstandard form. One can convert this to a standard form of \PIII with the 
transformation
\begin{equation}
\begin{split}
   t & = z^2
     \\
   u(t) & = v^2(z)\ ,
\end{split}
\end{equation}
and the parameters become $ \alpha = 0, \beta = 0, \gamma = -1, \delta = 1 $.
In the next step one requires that the \PIII system is transformed into our
original \PIIIprime system with the transformation
\begin{equation}
\begin{split}
   x & = z^2
     \\
   w(x) & = zv(z)\ ,
\end{split}
\label{IIIprime_xfm}
\end{equation}
and this does not alter the parameters. In the final step one has to restore the
original scaling of the \PIIIprime system with $ \gamma = 4, \delta = -4 $,
and all that is necessary is a scaling of the independent variable
$ T = \quarter\; x = \quarter\; t $. The other parameters $ \alpha, \beta $ are
unchanged and so we have our original system but with $ v_1 = -1, v_2 = 0 $.
The composition of all these transformations yields the stated result.
\end{proof}

\begin{proof}[Third Proof]
The inverse transformation can be written as
\begin{equation}
\begin{split}
   t & = 4T  
   \\
   q(t) & = -2TQ(T){ Q(T)-2TQ'(T) \over (Q^2(T)+T)(Q^2(T)-T) }\ ,
\end{split}
\label{Ixfm-riccati}
\end{equation}
and this can be recognised as a degenerative form of the Riccati transformations
\begin{equation}
\begin{split}
      z & = \mu t
   \\
   u(z) & = {v'(t)^2 + A_2(v)v'(t) + A_4(v) \over B_2(v)v'(t) + B_4(v)}
\end{split}
\end{equation}
first proposed by Fokas and Ablowitz \cite{FA_1982}, but generalised to include
a scale change in the independent variable. Here $ A_2(v), B_2(v) $ are
quadratic polynomials of $ v(t) $ with $t$-dependent coefficients and 
$ A_4(v), B_4(v) $ are quartic polynomials, and in this particular form the
$ (v')^2 $ and $ B_2(v)v' $ terms are absent.  One can derive (\ref{Ixfm-riccati}),
as well as the Gambier transformation for \PaII, using the approach of \cite{FA_1982}.
\end{proof}

\begin{remark}
The appearance of branch points at $ p = 0, 1 $ in the transformation 
(\ref{PIII_Wxfm}) is a consequence of the fixed singularities at these points 
in the differential equation (\ref{PV_dual}) that is dual to the 
general \PIII equation. Similarly the transformation for \PaII (\ref{PII_Hxfm}) has 
a branch point at $ p = 0 $ which corresponds to the singularity in 
the differential equation for PXIV. 
The square-root nature of the branches naturally appears from the inversion
of the quadratic rational transformation given by the last member of 
(\ref{PIII_Wixfm}).
\end{remark}

\begin{remark}
The solution of the \PIIIprime differential equation $ q(t) $ with parameters 
$ (0,0) $ is related to another solution $ Q(T) $ of this equation with
parameters $ (-1,0) $ by the direct transformation
\begin{equation}
\begin{split}
   T & = \quarter t 
   \\
   Q(T) & = i {\sqrt{t} \over 2q(t)}
   \left( \sqrt{tq'(t)-q^2(t)-t \over 2} + \sqrt{tq'(t)+q^2(t)-t \over 2}
   \right)\ .
\end{split}
\end{equation}
This follows by eliminating $ p $ from (\ref{PIII_Wxfm}) using the Hamilton 
equations (\ref{PIII_Heom}).
\end{remark}

\begin{remark}
Using the relations (\ref{IIIprime_xfm}) one can recast the transformations between the 
origin and $ (-1,0) $ for the \PIIIprime system into equivalent transformations
for the \PIII system. If one employs a notation for the 
\PIII transcendent $ \mu(s) $
with parameters $ \alpha = 0, \beta = 4, \gamma = 4, \delta = -4 $ and 
$ M(S) $ with parameters $ \alpha = 0, \beta = 0, \gamma = 4, \delta = -4 $
one has the direct and inverse transformations
\begin{align}
   S & = \half s
   \\
   M(S) & = {i \over 2\sqrt{s}\mu(s)}
        \Big[ \sqrt{\mu(s)+s\mu'(s)+2s(\mu^2(s)-1)} 
   \\
   & \phantom{= {i \over 2\sqrt{s}\mu(s)} \Big[}
             + \sqrt{\mu(s)+s\mu'(s)-2s(\mu^2(s)+1)} \Big]
   \nonumber\\
   \mu(s) & = {M(S)M'(S) \over (M^2(S)+1)(M^2(S)-1)} .
\end{align}
\end{remark}

\begin{remark}
One might think that transformations from the origin to other nearest neighbours
along the axes in parameter space (recall Fig. 2) would be simply expressible. 
One can indeed find these by 
composing B\"acklund transformations with the above transformation and has the
following Corollary. In the following discussion we employ a renaming of the 
variables $ Q, P, H \mapsto Q_W, P_W, H_W $ and a notation using the four
primary compass directions $ S, E, N, W $.
\end{remark}

\begin{corollary}
The canonical variables of the \PIIIprime Hamiltonian system 
$ \{q,p;t,h\} $ with $ {\bf v}=(0,0) $ 
are related to the \PIIIprime Hamiltonian system 
$ \{Q_S, P_S; T, H_S \} $ with $ {\bf v}=(0,-1) $ by the direct 
canonical transformation with $ T = \quarter t $ and
\begin{equation}
\begin{split}
   Q_S & = i {\sqrt{t} \over 2} {1 \over \sqrt{p} + \sqrt{p-1}}
                {\half + q\sqrt{p}\sqrt{p-1} - i\sqrt{t}\sqrt{p-1} \over
                 \half - q\sqrt{p}\sqrt{p-1} + i\sqrt{t}\sqrt{p-1} }
   \\
   P_S & = -{i \over \sqrt{t}} \left( \sqrt{p} + \sqrt{p-1} \right)
            \left[ \half - q\sqrt{p}\sqrt{p-1} + i\sqrt{t}\sqrt{p} \right] \ ,
\end{split}
\label{PIII_Sxfm}
\end{equation}
to the \PIIIprime Hamiltonian system $ \{Q_E, P_E; T, H_E \} $ with 
$ {\bf v}=(1,0) $ by the direct canonical transformation
\begin{equation}
\begin{split}
   Q_E & = -i {\sqrt{t} \over 2} \left( \sqrt{p} + \sqrt{p-1} \right)
               \left\{ 1
                 - {1 \over \half + q\sqrt{p}\sqrt{p-1} - i\sqrt{t}\sqrt{p-1} }
               \right.
   \\
          &    \hskip5cm\left.
                 - {1 \over \half + q\sqrt{p}\sqrt{p-1} + i\sqrt{t}\sqrt{p} }
               \right\}
   \\
   P_E & = -{i \over \sqrt{t}} 
            { \half + q\sqrt{p}\sqrt{p-1} + i\sqrt{t}\sqrt{p} \over
              \sqrt{p} + \sqrt{p-1} } \ ,
\end{split}
\label{PIII_Exfm}
\end{equation}
and to the \PIIIprime Hamiltonian system $ \{Q_N, P_N; T, H_N \} $ with 
$ {\bf v}=(0,1) $ by the direct canonical transformation
\begin{equation}
\begin{split}
   Q_N & = -i {\sqrt{t} \over 2} {1 \over \sqrt{p} + \sqrt{p-1}}
                {\half + q\sqrt{p}\sqrt{p-1} + i\sqrt{t}\sqrt{p} \over
                 \half - q\sqrt{p}\sqrt{p-1} - i\sqrt{t}\sqrt{p} }
   \\
   P_N & = {i \over \sqrt{t}} 
            { \half - q\sqrt{p}\sqrt{p-1} - i\sqrt{t}\sqrt{p} \over
              \sqrt{p} - \sqrt{p-1} } \ .
\end{split}
\label{PIII_Nxfm}
\end{equation}
\end{corollary}

\begin{proof}
These are found by composing the transformation (\ref{PIII_Wxfm}) with the 
B\"acklund transformations (\ref{PIII_shifts}).
\end{proof}

These transformations attain a clearer form when one considers the relationship
between the Hamiltonians defined by these canonical transformations. To illustrate
this we consider the case of \PaII first. 
Denoting $ H_{0} := H_{\rm II}(q,p;t;v_1 = 0) $ and
$ H_{1/2} := H_{\rm II}(Q,P;T,v_1 = \half) $ we find that by
employing the transformations (\ref{PII_xfm}) directly, or by constructing the 
generating function for the canonical transformation, that
\begin{align}
   H_{0} & = -2^{-1/3}\left[ P^2-2(Q^2+T/2)P+T^2/4 \right] \cr
         & = -2^{-1/3}\left[ 2H_{1/2}+Q+T^2/4 \right] \ .
\end{align}
But we note
$ H_{-1/2} := H_{\rm II}(T,v_1 = -\half) $ is given by the B\"acklund transformation
\begin{equation}
   H_{-1/2} = T_1 H_{1/2} = H_{1/2}+Q\ ,
\end{equation}
so that we have the additive relation
\begin{equation}
   H_{0} + t^2/8 = -2^{-1/3}( H_{1/2}+H_{-1/2} )\ ,
\end{equation}
and a multiplicative relation for the corresponding $\tau$-functions
\begin{equation}
   \exp(t^3/24) \tau_{0}(t) = \tau_{-1/2}(T)\tau_{1/2}(T) \ .
\end{equation}
If one takes the definition of the $\phi$-factors or $ \tau$-cocycles \cite{No_2000}, 
\begin{equation}
    \tau(t) := \exp(-t^3/24) \phi(t)\ ,
\end{equation}
then the multiplicative relation is 
\begin{equation}
   \phi_{0}(t) = \phi_{-1/2}(T)\phi_{1/2}(T) \ .
\end{equation}
The above identity has arisen in random matrix theory as a result of reconciling
two different calculations \cite{FW_2001a} for the cumulative distribution of the 
largest eigenvalue of a random hermitian matrix drawn from the Gaussian unitary 
ensemble.

We also have analogous relations in the case of \PIIIprime.

\begin{proposition}
The Hamiltonians for the parameters sets $ (-1,0) $, $ (0,-1) $, $ (1,0) $, 
$ (0,1) $ are related to the Hamiltonian at the origin by the additive relation
\begin{equation}
   TH_{W}(T)+TH_{S}(T)+TH_{E}(T)+TH_{N}(T) = th(t)-\quarter \ .
\label{PIII_Hsum}
\end{equation}
The canonical momenta are related in a similar way,
\begin{equation}
   P_{W}+P_{S}+P_{E}+P_{N} = 4p \ .
\label{PIII_Psum}
\end{equation}
The $\tau$-functions are related in the multiplicative form
\begin{equation}
   T^{1/4}\tau_{W}(T)\tau_{S}(T)\tau_{E}(T)\tau_{N}(T) = \tau(t)\ .
\label{PIII_tauproduct}
\end{equation}
\end{proposition}

\begin{proof}
By employing the canonical transformations 
(\ref{PIII_Wxfm},\ref{PIII_Sxfm},\ref{PIII_Exfm},\ref{PIII_Nxfm})
we find the Hamiltonians at the four corners $ (-1,0) $, $ (0,-1) $, $ (1,0) $, 
$ (0,1) $ are related to the Hamiltonian at the origin by
\begin{equation}
\begin{split}
  TH_{W}(T) & = 
  \quarter\; th(t)-\sixteenth+\quarter
           \left[ -\half-q\sqrt{p}\sqrt{p-1}+i\sqrt{t}(\sqrt{p}-\sqrt{p-1}) \right]
  \cr
  TH_{S}(T) & = 
  \quarter\; th(t)-\sixteenth+\quarter
           \left[  \half+q\sqrt{p}\sqrt{p-1}-i\sqrt{t}(\sqrt{p}+\sqrt{p-1}) \right]
  \cr
  TH_{E}(T) & = 
  \quarter\; th(t)-\sixteenth+\quarter
           \left[ -\half-q\sqrt{p}\sqrt{p-1}-i\sqrt{t}(\sqrt{p}-\sqrt{p-1}) \right]
  \cr
  TH_{N}(T) & = 
  \quarter\; th(t)-\sixteenth+\quarter
           \left[  \half+q\sqrt{p}\sqrt{p-1}+i\sqrt{t}(\sqrt{p}+\sqrt{p-1}) \right]\ ,
\end{split}
\end{equation}
and their sum yields (\ref{PIII_Hsum}).
Utilising the second members of 
(\ref{PIII_Wxfm},\ref{PIII_Sxfm},\ref{PIII_Exfm},\ref{PIII_Nxfm}) relation
(\ref{PIII_Psum}) follows directly. The final product relation for the
$\tau$-functions is the result of integrating (\ref{PIII_Hsum}).
\end{proof}

\medskip
\noindent{\it Acknowledgements}\\
It is a pleasure to acknowledge the many and wide-ranging discussions with
Peter Forrester, the valuable advice and insight of Chris Cosgrove and the 
observations made by Nalini Joshi and Peter Clarkson.
This research has been supported by the Australian Research Council.

\bibliographystyle{amsplain}
\bibliography{nonlinear,random_matrices}

\end{document}